\documentclass[reqno]{amsart}

\usepackage{amsmath}
\usepackage{amssymb}
\usepackage{amsthm}
\usepackage{epsfig}

\newtheorem{theorem}{Theorem}[section]
\newtheorem{prop}[theorem]{Proposition}
\newtheorem{lem}[theorem]{Lemma}
\newtheorem*{cor}{Corollary}
\theoremstyle{definition}
\newtheorem{defn}[theorem]{Definition}
\newtheorem*{notation}{Notation}
\theoremstyle{remark}
\newtheorem*{rem}{Remark}
\numberwithin{equation}{section}

\newcommand{\Gg}{\mathfrak{g}}    
\newcommand{\Gh}{\mathfrak{h}}

\newcommand{\Gr}{\mathfrak{r}}
\newcommand{\Gz}{\mathfrak{z}}

\begin{document}

\title[Quantum Heisenberg group algebra]
{Dressing orbits and a quantum Heisenberg group algebra}

\author{Byung-Jay Kahng}
\date{}
\address{Department of Mathematics\\ University of Kansas\\
Lawrence, KS 66045}
\email{bjkahng@math.ku.edu}

\begin{abstract}
In this paper, as a generalization of Kirillov's orbit theory, we
explore the relationship between the dressing orbits and irreducible
${}^*$-representations of the Hopf $C^*$-algebras $(A,\Delta)$ and
$(\tilde{A},\tilde{\Delta})$ we constructed earlier.  We discuss
the one-to-one correspondence between them, including their
topological aspects.

On each dressing orbit (which are symplectic leaves of the underlying
Poisson structure), one can define a Moyal-type deformed product at
the function level.  The deformation is more or less modeled by the 
irreducible representation corresponding to the orbit.  We point out
that the problem of finding a direct integral decomposition of the 
regular representation into irreducibles (Plancherel theorem) has an 
interesting interpretation in terms of these deformed products.
\end{abstract}
\maketitle

{\sc Introduction.}
According to Kirillov's orbit theory \cite{Ki2}, \cite{Ki}, the
representation theory of a Lie group is closely related with the
coadjoint action of the group on the dual vector space of its Lie
algebra.  The coadjoint orbits play the central role.  The program
is most successful for nilpotent or exponential solvable Lie groups.
While it does not work as well for other types of Lie groups, the
orbit theory, with some modifications, is still a very useful tool
in the Lie group representation theory.

 It is reasonable to expect that a generalization to some extent of 
the Kirillov-type orbit theory will exist even for quantum groups. 
On the other hand, it has been known for some time that the ``geometric 
quantization'' of physical systems is very much related with the 
construction of irreducible unitary representations in mathematics
(See \cite{Ko}, \cite{BFFL}, \cite{Ki4}.).  The orbit theoretical
approach is instrumental in these discussions.

 Generalizing the orbit theory to the quantum group level is a quite
interesting program \cite{Ki5}, and is still on-going since late 80's
(\cite{SV}, \cite{LS}).  In this paper, we will focus our attention 
to the examples of non-compact quantum groups $(A,\Delta)$ and $(\tilde{A},
\tilde{\Delta})$, which we constructed in \cite{BJKp2} (See also \cite
{BJKhj}, \cite{BJKppha}.).  Noting the fact that the classical counterparts
to $(A,\Delta)$ and $(\tilde{A},\tilde{\Delta})$ are exponential solvable
Lie groups, we plan to study these examples from the point of view of the 
(generalized) orbit theory.  Most of the results in this paper are not 
necessarily surprising.  But we still believe this is a worthwhile project, 
especially since we later plan to explore the aspect of our examples as 
certain quantized spaces as well as being quantum groups.

 Our goal here is twofold.  One is to study the representation theory 
of the examples $(A,\Delta)$ and $(\tilde{A},\tilde{\Delta})$, by using 
the orbit theory.  Although we now work with dressing actions instead 
of coadjoint actions, we will see that many of the results are analogous 
to those of the classical counterparts.  More long-term goal is to further 
investigate how the orbits and representation theory are related with the 
quantization process.

 This paper is organized as follows.  In section 1, we review some
preliminary results from our earlier papers.  We discuss the Poisson
structures, dressing actions and dressing orbits.  Then we recall the
definitions of our main examples $(A,\Delta)$ and $(\tilde{A},\tilde
{\Delta})$.  In section 2, we study the irreducible ${}^*$-representations
of our examples, in connection with the dressing orbits.  Going one step
further than just pointing out the one-to-one correspondence between
them, we will discuss the correspondence in terms of the topological
structures on the set of orbits and on the set of representations.

 For each dressing orbit ${\mathcal O}$, which is actually a symplectic 
leaf, we can construct a Moyal-type, deformed (quantized) product on
the space of smooth functions on ${\mathcal O}$.  We do this by
considering these smooth functions as operators on a Hilbert space.
We first find, in section~3, a canonical measure on each orbit, which
plays an important role in the construction of the Hilbert space.
The orbit deformation is carried out in section~4.  It turns out that
the deformation of an orbit ${\mathcal O}$ is ``modeled'' by the
irreducible ${}^*$-representation corresponding to ${\mathcal O}$.
Furthermore, we will see that the ``regular representation'' $L$,
which we used in \cite{BJKp2} to give the specific operator algebra
realization of $A$, is equivalent to a direct integral of the irreducible
${}^*$-representations.  This is a version of the Plancherel theorem.

 The discussion in section 4 (concerning the deformation of orbits) is
restricted to the case of $(A,\Delta)$, while the case of $(\tilde{A},
\tilde{\Delta})$ is postponed to a future paper.  However, we include
an Appendix at the end of section 4, where we give a short preliminary
report on the case of $(\tilde{A},\tilde{\Delta})$.  Although most of
the general ideas do go through, there are some technical obstacles
which we have to consider.  We hope that the discussions in sections 3
and 4 (as well as the ones in Appendix) will be helpful in our attempts
to understand the relationship between the orbits and the quantization
process.

\section{Preliminaries}

 Our objects of study are the Poisson--Lie groups $H$, $G$ and $\tilde{H}$,
$\tilde{G}$, as well as their quantizations (i.\,e. ``quantum groups'')
$(A,\Delta)$ and $(\tilde{A},\tilde{\Delta})$.  For more complete
descriptions of these objects, see \cite{BJKp2} (and also \cite{BJKhj}).
Let us begin with a short summary.

\subsection{The Poisson--Lie groups $H$, $G$, $\tilde{H}$, $\tilde{G}$.
The dressing orbits.}
The group $H$ is the $(2n+1)$-dimensional Heisenberg Lie group.  Its
underlying space is $\mathbb{R}^{2n+1}$ and the multiplication on it
is defined by
$$
(x,y,z)(x',y',z')=(x+x',y+y',z+z'+x\cdot y'),
$$
for $x,x',y,y'\in\mathbb{R}^n$ and $z,z'\in\mathbb{R}$.  We also consider
the extended Heisenberg group $\tilde{H}$, whose group law is given by
$$
(x,y,z,w)(x',y',z',w')=(x+e^w x',y+e^{-w} y',z+z'+e^{-w}x\cdot y',w+w').
$$
It is $(2n+2)$-dimensional.  The notation is similar as above, with
$w,w'\in\mathbb{R}$.  It is easy to see that $\tilde{H}$ contains $H$
as a normal subgroup.

 In \cite{BJKp2}, we obtained the ``dual Poisson--Lie group'' $\tilde{G}$
of $\tilde{H}$.  It is $(2n+2)$-dimensional, considered as a dual vector
space of $\tilde{H}$, and is determined by the multiplication law:
$$
(p,q,r,s)(p',q',r',s')=(e^{\lambda r'}p+p',e^{\lambda r'}q+q',r+r',s+s').
$$
And the dual Poisson--Lie group $G$ of $H$ is determined by the
multiplication law:
$$
(p,q,r)(p',q',r')=(e^{\lambda r'}p+p',e^{\lambda r'}q+q',r+r').
$$

\begin{rem}
In the above, $\lambda\in\mathbb{R}$ is a fixed constant, which determines
a certain non-linear Poisson structure when $\lambda\ne0$.  In section 1 of
\cite{BJKp2}, we gave a discussion on how the above pairs of Poisson--Lie
groups are related with a so-called ``classical $r$-matrix'' element.
Meanwhile, note that by taking advantage of the fact that the groups are
exponential solvable, we are considering them as vector spaces (identified
with the corresponding Lie algebras).
\end{rem}

 Given a dual pair of Poisson--Lie groups, there exists the so-called
{\em dressing action\/} of a Poisson--Lie group acting on its dual
Poisson--Lie group.  It is rather well known that the notion of a
dressing action is the natural generalization of the coadjoint action
of a Lie group acting on the dual space of its Lie algebra.  Also by
a result of Semenov-Tian-Shansky, it is known that the dressing orbits
are exactly the symplectic leaves in the Poisson--Lie groups \cite{LW},
\cite{Se}.

 It is customary to define the dressing action as a right action.  But
for the purpose of this paper and the future projects in our plans, it
is actually more convenient to work with the ``left'' dressing action.
It is related to our specific choice in \cite{BJKp2} of the multiplications
on $G$ and $(A,\Delta)$, so that the left Haar measure naturally comes
from the ordinary Lebesgue measure on $G$ (See also \cite{BJKppha}.).
To compute the left dressing action of $H$ on $G$ (similarly, the action
of $\tilde{H}$ on $\tilde{G}$), it is useful to consider the following
``double Lie group'' $\tilde{H}\Join\tilde{G}$.  It is isomorphic to
the definition considered in \cite{BJKhj}.

 Here and throughout this paper, we denote by $\eta_{\lambda}(r)$ the
expression, $\eta_{\lambda}(r):=\dfrac{e^{2\lambda r}-1}{2\lambda}$.
When $\lambda=0$, we take $\eta_{\lambda}(r)=r$.

\begin{lem}
Let $\tilde{H}\Join\tilde{G}$ be defined by the following multiplication
law:
\begin{align}
&(x,y,z,w;p,q,r,s)(x',y',z',w';p',q',r',s')  \notag \\
&=\bigl(x+e^{\lambda r+w}x',y+e^{\lambda r-w}y',  \notag \\
&\qquad z+z'+e^{\lambda r-w}x\cdot y'-\lambda p\cdot x'-\lambda q\cdot y'
+\lambda\eta_{\lambda}(r)x'\cdot y',w+w';   \notag \\
&\qquad e^{\lambda r'+w'}p+p'-e^{\lambda r'+w'}\eta_{\lambda}(r)y',
e^{\lambda r'-w'}q+q'+e^{\lambda r'-w'}\eta_{\lambda}(r)x', \notag \\
&\qquad r+r',s+s'-p\cdot x'+q\cdot y'+\eta_{\lambda}(r)x'\cdot y'\bigr).
\notag
\end{align}
We recover the group structures of $\tilde{H}$ and $\tilde{G}$,
by identifying $(x,y,z,w)\in\tilde{H}$ with $(x,y,z,w;0,0,0,0)$ and
$(p,q,r,s)\in\tilde{G}$ with $(0,0,0,0;p,q,r,s)$.  It is clear that
$\tilde{H}$ and $\tilde{G}$ are closed Lie subgroups.  Note also that
any element $(x,y,z,w;p,q,r,s)\in\tilde{H}\Join\tilde{G}$ can be
written as
$$
(x,y,z,w;p,q,r,s)=(x,y,z,w;0,0,0,0)(0,0,0,0;p,q,r,s).
$$
In other words, $\tilde{H}\Join\tilde{G}$ is the ``double Lie group''
of $\tilde{H}$ and $\tilde{G}$.  Meanwhile, if we consider only the
$(x,y,z)$ and the $(p,q,r)$ variables, we obtain in the same way the
double Lie group $H\Join G$ of $H$ and $G$.
\end{lem}

\begin{proof}
The group is obtained by first considering the ``double Lie algebra''
(\cite{LW}), which comes from the Lie bialgebra $(\tilde{\Gh},
\tilde{\Gg})$ corresponding to the Poisson structures on $\tilde{H}$
and $\tilde{G}$ (See \cite{BJKhj} for computation.).  Verification
of the statements are straightforward.
\end{proof}

 The left dressing action is defined exactly in the same manner as
the usual (right) dressing action.  For $h\in\tilde{H}$ and $\mu
\in\tilde{G}$, regarded naturally as elements in $\tilde{H}\Join
\tilde{G}$, we first consider the product $\mu\cdot h$.  Factorize
the product as $\mu\cdot h=h^{\mu}\cdot\mu^{h}$, where $h^{\mu}
\in\tilde{H}$ and $\mu^{h}\in\tilde{G}$.  The left dressing action,
$\delta$, of $\tilde{H}$ on $\tilde{G}$ is then given by $\delta_h
(\mu):=\mu^{(h^{-1})}$.

 We are mainly interested in the dressing orbits contained in $G$
and $\tilde{G}$, which will be the symplectic leaves.  The following
two propositions give a brief summary.  For the computation of the
Poisson bracket on $G$, see Theorem 2.2 of \cite{BJKp2}.  The Poisson
bracket on $\tilde{G}$ is not explicitly mentioned there, but we can
more or less follow the proof for the case of $G$.  The computations
for the dressing orbits are straightforward from the definition of
the Poisson brackets.

\begin{prop}\label{dressingG}
\begin{enumerate}
\item The Poisson bracket on $G$ is given by the following expression:
$$
\{\phi,\psi\}(p,q,r)=\eta_{\lambda}(r)(x\cdot y'-x'\cdot y),
\qquad{\text {for\ }}\phi,\psi\in C^{\infty}(G).
$$
Here $d\phi(p,q,r)=(x,y,z)$ and $d\psi(p,q,r)=(x',y',z')$, which are
naturally considered as elements of $\Gh$.
\item The left dressing action of $H$ on $G$ is:
$$
\bigl(\delta(a,b,c)\bigr)(p,q,r)=\bigl(p+\eta_{\lambda}(r)b,
q-\eta_{\lambda}(r)a,r\bigr).
$$
\item The dressing orbits in $G$ are:
\begin{itemize}
\item ${\mathcal O}_{p,q}=\bigl\{(p,q,0)\bigr\}$, when $r=0$.
\item ${\mathcal O}_r=\bigl\{(\alpha,\beta,r):(\alpha,\beta)\in
\mathbb{R}^{2n}\bigr\}$, when $r\ne0$.
\end{itemize}
The ${\mathcal O}_{p,q}$ are 1-point orbits and the ${\mathcal O}_r$ are
$2n$-dimensional orbits.  By standard theory, these orbits are exactly
the symplectic leaves in $G$ for the Poisson bracket on it defined above.
\end{enumerate}
\end{prop}

\begin{prop}\label{dressingGG}
\begin{enumerate}
\item The Poisson bracket on $\tilde{G}$ is given by
$$
\{\phi,\psi\}(p,q,r,s)=p\cdot(wx'-w'x)+q\cdot(w'y-wy')+\eta_{\lambda}(r)
(x\cdot y'-x'\cdot y),
$$
for $\phi,\psi\in C^{\infty}(\tilde{G})$.  We are again using the
natural identification of $d\phi(p,q,r,s)=(x,y,z,w)$ and $d\psi(p,q,r,s)
=(x',y',z',w')$ as elements of $\tilde{\Gh}$.
\item The left dressing action of $\tilde{H}$ on $\tilde{G}$ is:
\begin{align}
\bigl(\delta(a,b,c,d)\bigr)(p,q,r,s)&=\bigl(e^{-d}p+\eta_{\lambda}(r)b,
e^{d}q-\eta_{\lambda}(r)a,r,  \notag \\
&\qquad s+e^{-d}p\cdot a-e^{d}q\cdot b+\eta_{\lambda}(r)a\cdot b\bigr).  \notag
\end{align}
\item The dressing orbits in $\tilde{G}$ are:
\begin{itemize}
\item $\tilde{\mathcal O}_s=\bigl\{(0,0,0,s)\bigr\}$, when $(p,q,r)=(0,0,0)$.
\item $\tilde{\mathcal O}_{p,q}=\bigl\{(\alpha p,\frac{1}{\alpha}q,0,\gamma):
\alpha>0,\gamma\in\mathbb{R}\bigr\}$, when $r=0$ but $(p,q)\ne(0,0)$.
\item $\tilde{\mathcal O}_{r,s}=\bigl\{(\alpha,\beta,r,s-\frac{1}
{\eta_{\lambda}(r)}\alpha\cdot\beta):(\alpha,\beta)\in\mathbb{R}^{2n}\bigr\}$,
when $r\ne0$.
\end{itemize}
The $\tilde{\mathcal O}_s$ are 1-point orbits, the $\tilde{\mathcal O}
_{p,q}$ are 2-dimensional orbits, and the $\tilde{\mathcal O}_{r,s}$
are $2n$-dimensional orbits.  These are exactly the symplectic leaves in
$\tilde{G}$ for the Poisson bracket on it.
\end{enumerate}
\end{prop}

\begin{rem}
We should point out that the notation for the 2-dimensional orbits $\tilde
{\mathcal O}_{p,q}$ in $\tilde{G}$ are somewhat misleading.  We can see
clearly that we can have $\tilde{\mathcal O}_{p,q}=\tilde{\mathcal O}_{p',q'}$,
if $p'=\alpha p$ and $q'=\frac{1}{\alpha}q$ for some $\alpha>0$.  To avoid
introducing cumbersome notations, we nevertheless chose to stay with this
ambiguity.  A sketch is given below (in figure 1) to help us visualize
the situation.

\center{\epsfig{file=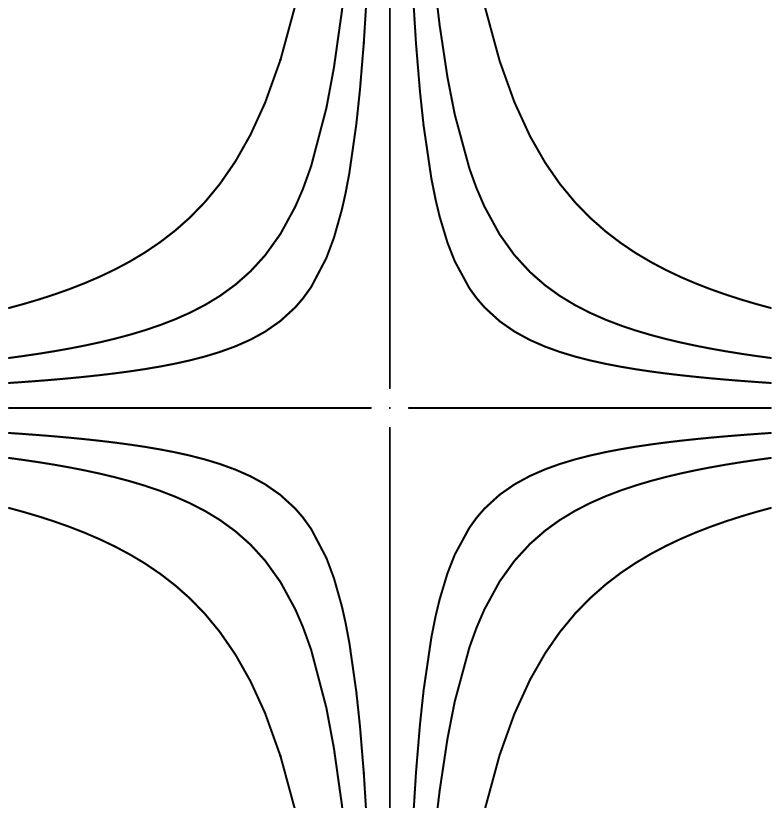,width=1.5in}}

(Figure 1: Trace curves of the orbits on the $r=s=0$ plane)
\end{rem}

\subsection{The Hopf $C^*$-algebras $(A,\Delta)$ and $(\tilde{A},
\tilde{\Delta})$.}
By realizing that the Poisson structures on $G$ is a non-linear Poisson bracket
of the ``cocycle perturbation'' type (as in \cite{BJKp1}), we were able to
construct the Hopf $C^*$-algebras (quantum groups) $(A,\Delta)$ and $(\tilde{A},
\tilde{\Delta})$ by deformation.  For their precise definitions, see \cite{BJKp2}.

 As a $C^*$-algebra, $A$ is isomorphic to a twisted crossed product algebra.
That is, $A\cong C^*\bigl(H/Z,C_0(G/{Z^{\perp}}),\sigma\bigr)$, where $H$ and
$G(=H^*)$ are as in \S1.1 and $Z$ is the center of $H$ (so $Z=\bigl\{{\text
{$(0,0,z)$'s}}\bigr\}$).  We denoted by $\sigma$ the twisting cocycle for
the group $H/Z$.  As constructed in \cite{BJKp2}, $\sigma$ is a continuous
field of cocycles $G/{Z^{\perp}}\ni r\mapsto\sigma^r$, where
$$
\sigma^r\bigl((x,y),(x',y')\bigr)=\bar{e}\bigl[\eta_{\lambda}(r)x\cdot y'\bigr].
$$
Following the notation of the previous papers, we are letting $e(t)=
e^{(2\pi i)t}$ and $\bar{e}(t)=e^{(-2\pi i)t}$, while $\eta_{\lambda}(r)$
is as before.  The elements $(x,y)$, $(x',y')$ are group elements in $H/Z$.
Via a certain ``regular representation'' $L$, we were able to realize the
$C^*$-algebra $A$ as an operator algebra in ${\mathcal B}({\mathcal H})$,
where ${\mathcal H}=L^2(H/Z\times G/{Z^{\perp}})$ is the Hilbert space
consisting of the $L^2$-functions in the $(x,y,r)$ variables.

 In \cite{BJKp2}, we showed that the $C^*$-algebra $A$ is a strict deformation
quantization (in the sense of Rieffel \cite{Rf4}) of $C_0(G)$.  For convenience,
the deformation parameter $\hbar$ has been fixed ($\hbar=1$), which is the
reason why we do not see it in the definition of $A$.  When $\hbar=0$ (i.\,e.
classical limit), we take $\sigma\equiv1$.  Then $A_{\hbar=0}\cong C_0(G)$.
Throughout this paper (as in our previous papers), we write $A=A_{\hbar=1}$.
On $A$, an appropriate comultiplication can be defined using a certain ``(regular)
multiplicative unitary operator''.  Actually, we can show (see \cite{BJKppha})
that $(A,\Delta)$ is an example of a {\em locally compact quantum group\/},
in the sense of Kustermans and Vaes \cite{KuVa}.  All these can be done
similarly for the case of $\tilde{H}$ and $\tilde{G}$, obtaining $(\tilde{A},
\tilde{\Delta})$.

 Since our goal here is to study the ${}^*$-representations of $A$ and
$\tilde{A}$, we will go lightly on the discussion of their quantum group
structures.  We will review the notations as the needs arise.  For the
most part, it will be useful to recall that we can regard $(A,\Delta)$ as
a ``quantized $C^*(H)$'', i.\,e. a ``quantum Heisenberg group algebra''.
Or dually, we may regard it as a ``quantized $C_0(G)$''.  Similar comments
hold for $(\tilde{A},\tilde{\Delta})$, which can be considered as an ``extended
quantum Heisenberg group algebra'' (``quantized $C^*(\tilde{H})$'') or as
a ``quantized $C_0(\tilde{G})$''.

\section{The irreducible ${}^*$-representations}

 The irreducible ${}^*$-representations of $(A,\Delta)$ and $(\tilde{A},
\tilde{\Delta})$ have been found in \cite{BJKhj}, taking advantage of the
fact that $A$ and $\tilde{A}$ are twisted group $C^*$-algebras.  The results
are summarized in the first two propositions below.  Note that we study here
the ${}^*$-representations of $(A,\Delta)$ and $(\tilde{A},\tilde{\Delta})$,
instead of their corepresentations. By the observation given in the previous
section, the ${}^*$-representations of $A$ [and $\tilde{A}$] more or less
correspond to the group representations of $H$ [and $\tilde{H}$].

\begin{notation}
In the below, ${\mathcal A}$ is the space of Schwartz functions in the
$(x,y,r)$ variables, having compact support in the $r$ variable.  Similarly,
$\tilde{\mathcal A}$ is the space of Schwartz functions in the $(x,y,r,w)$
variables, having compact support in the $r$ and $w$.  By viewing these
functions as operators (via the regular representation $L$; see equation
(2.4) and Example 3.6 of \cite{BJKp2}), we saw in the previous papers that
they are dense ${}^*$-subalgebras of $A$ and $\tilde{A}$.  Most of our specific
calculations have been carried out at the level of these function algebras.
\end{notation}

\begin{prop}\label{repA}
Every irreducible representation of ${\mathcal A}$ is equivalent to one of
the following representations.  Here, $f\in{\mathcal A}$.
\begin{itemize}
\item For $(p,q)\in\mathbb{R}^{2n}$, there is a 1-dimensional representation
$\pi_{p,q}$ of ${\mathcal A}$, defined by
$$
\pi_{p,q}(f)=\int f(x,y,0)\bar{e}(p\cdot x+q\cdot y)\,dxdy.
$$
\item For $r\in\mathbb{R}$, $r\ne0$, there is a representation $\pi_r$ of
${\mathcal A}$, acting on the Hilbert space ${\mathcal H}_r=L^2(\mathbb{R}^n)$
and is defined by
$$
\pi_r(f)\xi(u)=\int f(x,y,r)\bar{e}\bigl(\eta_{\lambda}(r)u\cdot y\bigr)
\xi(u+x)\,dxdy.
$$
\end{itemize}
We thus obtain all the irreducible ${}^*$-representations of $A$ by naturally
extending these representations.  We will use the same notation, $\pi_{p,q}$
and $\pi_r$, for the representations of $A$ constructed in this way.
\end{prop}

\begin{prop}\label{repAA}
The irreducible ${}^*$-representations of $\tilde{A}$ are obtained by
naturally extending the following irreducible representations of the
dense subalgebra $\tilde{\mathcal A}$.  Here, $f\in\tilde{\mathcal A}$.
\begin{itemize}
\item For $s\in\mathbb{R}$, there is a 1-dimensional representation
$\tilde{\pi}_s$ defined by
$$
\tilde{\pi}_s(f)=\int f(x,y,0,w)\bar{e}(sw)\,dxdydw.
$$
\item For $(p,q)\in\mathbb{R}^{2n}$, $(p,q)\ne(0,0)$, there is a representation
$\tilde{\pi}_{p,q}$ acting on the Hilbert space $\tilde{\mathcal H}_{p,q}=L^2
(\mathbb{R})$ defined by
$$
\tilde{\pi}_{p,q}(f)\zeta(d)=\int f(x,y,0,w)\bar{e}(e^{d}p\cdot x
+e^{-d}q\cdot y)\zeta(d+w)\,dxdydw.
$$
\item For $(r,s)\in\mathbb{R}^2$, $r\ne0$, there is a representation
$\tilde{\pi}_{r,s}$ acting on the Hilbert space $\tilde{\mathcal H}_{r,s}
=L^2(\mathbb{R}^n)$ defined by
$$
\tilde{\pi}_{r,s}(f)\xi(u)=\int f(x,y,r,w)\bar{e}(sw)\bar{e}\bigl
(\eta_{\lambda}(r)u\cdot y\bigr)(e^{-\frac{w}{2}})^n\xi(e^{-w}u
+e^{-w}x)\,dxdydw.
$$
\end{itemize}
We will use the same notation, $\tilde{\pi}_s$, $\tilde{\pi}_{p,q}$,
and $\tilde{\pi}_{r,s}$, for the corresponding representations of
$\tilde{A}$.
\end{prop}

\begin{rem}
For the construction of these representations, see section 2 of \cite{BJKhj}.
Meanwhile, a comment similar to an earlier remark has to be made about the
representations $\tilde{\pi}_{p,q}$.  It is not very difficult to see that
$\tilde{\pi}_{p,q}$ and $\tilde{\pi}_{p',q'}$ are equivalent if and only if
$p'=e^{l}p$ and $q'=e^{-l}q$ for some real number $l$.  Again, to avoid
having to introduce cumbersome notations, we are staying with the (possibly
ambiguous) notation used above.
\end{rem}

 These are all the irreducible ${}^*$-representations of $A$ and $\tilde{A}$
up to equivalence (assuming we accept the ambiguity mentioned in the above
remark).  We did not rely on the dressing orbits to find these representations
(we constructed the irreducible representations via the machinery of induced
representations \cite{BJKhj}), but we can still observe that the irreducible
representations of $A$ and $\tilde{A}$ are in one-to-one correspondence with
the dressing orbits in $G$ and $\tilde{G}$, respectively.  To emphasize the
correspondence, we used the same subscripts for the orbits and the related
irreducible representations.

 Let us denote by ${\mathcal O}(A)$ the set of dressing orbits contained
in $G$.  Since $G$ is being identified with its Lie algebra $\Gg$, it is
equipped with the $(2n+1)$-dimensional vector space topology.  Note that
on $G$, we have an equivalence relation such that $(p,q,r)\sim(p',q',r')$
if they are contained in the same orbit.  By viewing ${\mathcal O}(A)=
G/{\sim}$, we can give ${\mathcal O}(A)$ a natural quotient topology.

 Meanwhile, let $\operatorname{Irr}(A)$ be the set of equivalence classes
of irreducible ${}^*$-representations of $A$ (Proposition \ref{repA}).
For every representation $\pi\in\operatorname{Irr}(A)$, its kernel is
a primitive ideal of $A$.  Consider the ``Jacobson topology'' \cite{D1}
on $\operatorname{Prim}(A)$, that is, the closure of a subset ${\mathcal U}
\subseteq\operatorname{Prim}(A)$ is defined to be the set of all ideals
in $\operatorname{Prim}(A)$ containing the intersection of the elements
of ${\mathcal U}$.  Since the map $\pi\mapsto\operatorname{Ker}(\pi)$ is
a canonical surjective map from $\operatorname{Irr}(A)$ onto $\operatorname
{Prim}(A)$, the Jacobson topology on $\operatorname{Prim}(A)$ can be
carried over to $\operatorname{Irr}(A)$ (In our case, we will actually
have $\operatorname{Irr}(A)\cong\operatorname{Prim}(A)$, since $A$ is
a type I $C^*$-algebra.).

 In exactly the same way, we can also define ${\mathcal O}(\tilde{A})$
and $\operatorname{Irr}(\tilde{A})$ together with their respective
topological structures.  We already know that ${\mathcal O}(A)\cong
\operatorname{Irr}(A)$ and ${\mathcal O}(\tilde{A})\cong\operatorname
{Irr}(\tilde{A})$ as sets.  Let us now explore these correspondences
a little further, in terms of their respective topologies.

 For ordinary Lie groups, the question of establishing a topological
homeomorphism between the (coadjoint) orbit space ${\mathcal O}(G)$
and the representation space $\operatorname{Irr}\bigl(C^*(G)\bigr)$
is called by some authors as the ``Kirillov conjecture''.  It is
certainly well known to be true in the case of nilpotent Lie groups.
It is also true in the case of exponential solvable Lie groups (The
proof was established rather recently \cite{LL}.).

 In our case, the objects of our study are Hopf $C^*$-algebras (quantum
groups), and we consider dressing orbits instead of coadjoint orbits.
On the other hand, their classical limits are exponential solvable Lie
groups.  Because of this, for the proof of $\operatorname{Irr}(\tilde{A})
\cong{\mathcal O}(\tilde{A})$ and $\operatorname{Irr}(A)\cong{\mathcal O}
(A)$, it is possible to take advantage of the general result at the Lie
group setting.  Here is the main result, which is not really surprising:

\begin{theorem}\label{orbitrepcorr}
\begin{enumerate}
\item $\operatorname{Irr}(\tilde{A})\cong{\mathcal O}(\tilde{A})$, as
topological spaces.
\item $\operatorname{Irr}(A)\cong{\mathcal O}(A)$, as topological
spaces.
\end{enumerate}
\end{theorem}

\begin{proof}
If we look at the ${}^*$-representations in $\operatorname{Irr}(\tilde
{A})$ computed earlier, we can see that they closely resemble the
${}^*$-representations in $\operatorname{Irr}\bigl(C^*(\tilde{H})\bigr)$,
the equivalence classes of unitary group representations of the (exponential)
Lie group $\tilde{H}$.  It actually boils down to replacing ``$\eta_{\lambda}
(r)$'' with ``$r$''.  The two sets are certainly different, but noting that
$r\mapsto\eta_{\lambda}(r)\left(=\frac{e^{2\lambda r}-1}{2\lambda}\right)$
is one-to-one and onto, it is clear that the topological structures on
$\operatorname{Irr}(\tilde{A})$ and on $\operatorname{Irr}\bigl(C^*(\tilde{H})
\bigr)$ are the same.  Similar comment holds for the dressing orbit space
${\mathcal O}(\tilde{A})$ and the coadjoint orbit space ${\mathcal O}(\tilde{H})$.

By general theory on representations of exponential Lie groups (see section 3
of \cite{LL}), we have ${\mathcal O}(\tilde{H})\cong\operatorname{Irr}
\bigl(C^*(\tilde{H})\bigr)$.  Therefore, it follows from the observation
in the previous paragraph that ${\mathcal O}(\tilde{A})\cong\operatorname{Irr}
(\tilde{A})$.  The proof for the homeomorphism ${\mathcal O}(A)\cong
\operatorname{Irr}(A)$ can be carried out in exactly the same way, this
time considering the (nilpotent) Lie group $H$.
\end{proof}

 This is the general result we wanted to establish.  But one drawback of
the above proof is that it is rather difficult to see what is actually
going on.  To illustrate and for a possible future use, we collect in the
below a few specific results (with proofs given by direct computations)
showing the topological properties on $\operatorname{Irr}(\tilde{A})$.
We do not mention the case of $A$ here, but it would be obviously simpler.

 Let us begin by describing the quotient topology on ${\mathcal O}(\tilde{A})$.
Consider the ``points'' (i.\,e. orbits) in ${\mathcal O}(\tilde{A})$.  For
$r\ne0$, the topology on the set of the points $\tilde{\mathcal O}_{r,s}$
is the standard one, which is just the topology on the $(r,s)$ plane
excluding the $r$-axis.  When $r=0$ and $(p,q)=(0,0)$, in which case the
points consist of the orbits $\tilde{\mathcal O}_s$, the topology is exactly
the standard topology on a line (the $s$-axis).  It is non-standard in the
case when $r=0$ and $(p,q)\ne(0,0)$, where the points consist of the orbits
$\tilde{\mathcal O}_{p,q}$.  To visualize, the picture (figure~1) given at
the end of \S1.1 will be helpful here.

 Let us now turn our attention to $\operatorname{Irr}(\tilde{A})$.  By Theorem
\ref{orbitrepcorr}, we already know that the topology on it coincides with
the quotient topology on ${\mathcal O}(\tilde{A})$, under the identification
of the two sets via our one-to-one correspondence.  In the following three
propositions, we give direct proofs of a few selected situations that show
some non-standard topological behavior, i.\,e. when $r=0$ and $(p,q)\ne(0,0)$.
The notation for the representations in $\operatorname{Irr}(\tilde{A})$ are
as before.

\begin{prop}\label{rep(pq00)}
Consider the sequence of representations $\{\tilde{\pi}_{p,q}\}$, letting
$(p,q)$ approach $(0,0)$.  Then the limit points of the sequence are the
representations $\tilde{\pi}_s$.
\end{prop}

\begin{proof}
As $(p,q)\to(0,0)$, the $\tilde{\pi}_{p,q}$ approach the (reducible)
representation $S$, which acts on the Hilbert space $L^2(\mathbb{R})$
and is defined by
$$
S(f)\zeta(d)=\int f(x,y,0,w)\zeta(d+w)\,dxdydw.
$$
To see how this representation $S$ decomposes into, consider the
unitary map (Fourier transform) on $L^2(\mathbb{R})$ given by
$$
{\mathcal F}\zeta(s)=\int\zeta(\tilde{d})\bar{e}(s\tilde{d})\,d\tilde{d}.
$$
Using ${\mathcal F}$, we can define the representation $\tilde{S}$
which is equivalent to $S$.  By a straightforward calculation
using Fourier inversion theorem, we have:
$$
\tilde{S}(f)\zeta(s)={\mathcal F}S(f){\mathcal F}^{-1}\zeta(s)=
\int f(x,y,0,w)e(sw)\zeta(s)\,dxdydw.
$$
We can see that $\tilde{S}$ is the direct integral of the irreducible
representations $\tilde{\pi}_s$.  In other words,
$$
S\cong\tilde{S}=\int^{\oplus}\tilde{\pi}_s\,ds.
$$
Therefore, $\operatorname{Ker}S\subseteq\bigcap_s\operatorname{Ker}
\tilde{\pi}_s$.  It follows that all the $\tilde{\pi}_s$ are limit
points of the sequence $\{\tilde{\pi}_{p,q}\}_{(p,q)\to(0,0)}$,
under the topology on $\operatorname{Irr}(\tilde{A})$.  In view of
the result of Theorem \ref{orbitrepcorr} that $\operatorname{Irr}
(\tilde{A})\cong{\mathcal O}(\tilde{A})$, they will exhaust all the
limit points.
\end{proof}

\begin{prop}
Consider the sequence of representations $\{\tilde{\pi}_{p,q}\}$, letting
$q\to0$ while we hold $p$.  More specifically, consider the sequence
$\{\tilde{\pi}_{p,c_nq}\}$, where $p$ and $q$ are fixed and $\{c_n\}$
is a sequence of positive numbers approaching $0$.  Then the limit points
are the representations $\tilde{\pi}_{p,0}$ and $\tilde{\pi}_{0,q}$.
\end{prop}

\begin{proof}
It is clear that $\tilde{\pi}_{p,0}$ is a limit point of the sequence
$\{\tilde{\pi}_{p,c_nq}\}$.  Meanwhile, note that $\tilde{\pi}_{p',q'}$
is equivalent to $\tilde{\pi}_{e^{l}p',e^{-l}q'}$ (See the Remark following
Proposition \ref{repAA}.).   Let us choose a sequence of real numbers
$\{l_n\}$ such that $e^{-l_n}c_n=1$.  Then each $\tilde{\pi}_{p,c_nq}$
is equivalent to $\tilde{\pi}_{e^{l_n}p,q}$.  Since the $e^{l_n}=c_n$
obviously approach $0$, we conclude that $\tilde{\pi}_{0,q}$ is also
a limit point of the sequence $\{\tilde{\pi}_{p,c_nq}\}$.
\end{proof}

\begin{rem}
Similarly, we may consider the sequence $\{\tilde{\pi}_{c_np,q}\}$,
where $p$ and $q$ are fixed and $\{c_n\}$ is a sequence of positive
numbers approaching $0$.  In exactly the same way as above, we can
show that the limit points are the representations $\tilde{\pi}_{p,0}$
and $\tilde{\pi}_{0,q}$.  By modifying the proof a little, we can also
obtain various results of similar flavor.
\end{rem}

\begin{prop}
Consider the sequence of representations $\{\tilde{\pi}_{r,s}\}$, letting
$r\to0$ while $s$ is fixed.  Then all the representations $\tilde{\pi}_{p,q}$
and the $\tilde{\pi}_s$ are limit points of the sequence.
\end{prop}

\begin{proof}
As $r\to0$, the $\tilde{\pi}_{r,s}$ approach the following (reducible)
representation, $T$, acting on the Hilbert space $L^2(\mathbb{R}^n)$:
$$
T(f)\xi(u)=\int f(x,y,0,w)\bar{e}(sw)(e^{-\frac{w}{2}})^n
\xi(e^{-w}u+e^{-w}x)\,dxdydw.
$$
To see how $T$ decomposes into, consider the Fourier transform on
$L^2(\mathbb{R}^n)$.
$$
{\mathcal F}\xi(\alpha)=\int\xi(u)\bar{e}(\alpha\cdot u)\,du.
$$
As before, we can define the representation $\tilde{T}$ which is
equivalent to $T$.  By a straightforward calculation involving
Fourier inversion theorem, we have:
$$
\tilde{T}(f)\xi(\alpha)={\mathcal F}T(f){\mathcal F}^{-1}\xi(\alpha)
=\int f(x,y,0,w)e(\alpha\cdot x)\bar{e}(sw)(e^{\frac{w}{2}})^n
\xi(e^w\alpha)\,dxdydw.
$$
Suppose we expressed $\alpha$ as $\alpha=e^dp$, for some $p\in
{\mathbb{R}^n}$ and $d\in{\mathbb{R}}$.  Then it becomes:
$$
\tilde{T}(f)\xi(e^dp)=\int f(x,y,0,w)e(e^dp\cdot x)\bar{e}(sw)
(e^{\frac{w}{2}})^n\xi(e^{d+w}p)\,dxdydw.
$$
Fix $p(\ne0)$, and let us write $\xi_p(d)$ to be $\xi_p(d):=(e^{\frac
{d}{2}})^n\xi(e^dp)$.  Then we now have:
$$
\tilde{T}(f)\xi_p(d)=\int f(x,y,0,w)e(e^dp\cdot x)\bar{e}(sw)
\xi_p(d+w)\,dxdydw.
$$
This is really the expression for the inner tensor product representation,
$\tilde{\pi}_{-p,0}\boxtimes\tilde{\pi}_s$, which is equivalent to $\tilde
{\pi}_{-p,0}$ (the proof is by straightforward calculation, similar to the
one given in Proposition 4.4 of \cite{BJKhj}). As in the proof of Proposition
\ref{rep(pq00)}, we thus have:
$\operatorname{Ker}T\subseteq\bigcap_p\operatorname{Ker}
\tilde{\pi}_{p,0}$.  It follows that all the $\tilde{\pi}_{p,0}$ are
limit points of the sequence $\{\tilde{\pi}_{r,s}\}_{r\to0}$, under
the topology on $\operatorname{Irr}(\tilde{A})$.  By the result of
Proposition \ref{rep(pq00)}, it also follows that all the $\tilde{\pi}_s$
are limit points of the sequence.

Meanwhile, to look for more limit points of the sequence, let us define
the representations $Q_{r,s}$ as follows, which are equivalent to
the $\tilde{\pi}_{r,s}$.  We consider $\tilde{\pi}_{r,s}$ and consider
the Hilbert space $L^2(\mathbb{R}^n)$ on which $\tilde{\pi}_{r,s}(\tilde
{A})$ acts.  In $L^2(\mathbb{R}^n)$, define the unitary map ${\mathcal F}_r$
defined by
$$
{\mathcal F}_r\xi(v)=\int\xi(u)\bar{e}\bigl(\eta_{\lambda}(r)u\cdot v
\bigr)\bigl\|\eta_{\lambda}(r)\bigr\|^{\frac{1}{2}}\,du.
$$
This is again a kind of a Fourier transform, taking advantage of the
existence of the bilinear form $(u,v)\mapsto\eta_{\lambda}(r)u\cdot v$
in $\mathbb{R}^n$.  Its inverse is given by
$$
{\mathcal F}_r^{-1}\xi(u)=\int\xi(v)e\bigl(\eta_{\lambda}(r)u\cdot v
\bigr)\bigl\|\eta_{\lambda}(r)\bigr\|^{\frac{1}{2}}\,dv.
$$
We define the representation $Q_{r,s}$ by $Q_{r,s}(f)\xi:={\mathcal F}_r
\tilde{\pi}_{r,s}(f){\mathcal F}_r^{-1}\xi$.  We then have:
$$
Q_{r,s}(f)\xi(v)=\int f(x,y,r,w)\bar{e}(sw)e\bigl[\eta_{\lambda}(r)
(v+y)\cdot x\bigr](e^{\frac{w}{2}})^n\xi(e^wv+e^wy)\,dxdydw.
$$
Since the $Q_{r,s}$ are equivalent to the $\tilde{\pi}_{r,s}$,
we may now consider the sequence $\{Q_{r,s}\}_{r\to0}$.  We can
see right away that as $r\to0$, the sequence approaches the following
(reducible) representation, $Q$, acting on $L^2(\mathbb{R}^n)$:
$$
Q(f)\xi(v)=\int f(x,y,0,w)\bar{e}(sw)(e^{\frac{w}{2}})^n
\xi(e^wv+e^wy)\,dxdydw.
$$
We can follow exactly the same method we used in the case of the
representation $T$ to show that we now have: $\operatorname{Ker}Q
\subseteq\bigcap_q\operatorname{Ker}\tilde{\pi}_{0,q}$.  It follows
that all the $\tilde{\pi}_{0,q}$ are limit points of the sequence
$\{Q_{r,s}\}_{r\to0}$, or equivalently the sequence $\{\tilde
{\pi}_{r,s}\}_{r\to0}$.

In this way, we have shown so far that the representations $\tilde
{\pi}_{p,0}$, the $\tilde{\pi}_{0,q}$, as well as the $\tilde
{\pi}_s$ are limit points of the sequence.  It will be somewhat
cumbersome, but by choosing some suitable realizations of the
representations $\tilde{\pi}_{r,s}$, it is possible to show that
all the representations $\tilde{\pi}_{p,q}$ are also limit points.
\end{proof}

 Observe that the topological behaviors of $\operatorname{Irr}(\tilde
{A})$ manifested in the above three propositions are exactly those
of the quotient topology on ${\mathcal O}(\tilde{A})$, as was to be
expected from Theorem \ref{orbitrepcorr}.  On the other hand, it is
not quite sufficient to claim only from these types of propositions
that the one-to-one correspondences $\operatorname{Irr}(\tilde{A})\cong
{\mathcal O}(\tilde{A})$ and $\operatorname{Irr}(A)\cong{\mathcal O}
(A)$ are topological homeomorphisms.  For this reason, the actual
proof was given indirectly.

 These results give affirmation that there is a strong analogy between
our ``quantum'' case and the ``classical'' case of ordinary groups.
This is the underlying theme of this article.  On the other hand,
see \cite{BJKhj}, where we discuss an interesting ``quantum'' behavior
enjoyed by the ${}^*$-representations (e.\,g. the {\em quasitriangular\/}
property), due to the role played by the comultiplications on $(A,\Delta)$
and $(\tilde{A},\tilde{\Delta})$.

\section{Canonical measure on an orbit}

 We will work mostly with the dressing orbits in ${\mathcal O}(\tilde{A})$.
The case of orbits in ${\mathcal O}(A)$ will be simpler.

 Recall first the (left) dressing action, $\delta$, of $\tilde{H}$
on $\tilde{G}$, as defined in Proposition \ref{dressingGG}. By our
identification of $\tilde{G}$ with its Lie algebra $\tilde{\Gg}$, the
dressing action $\delta$ can be viewed as an action of $\tilde{H}$ on
$\tilde{\Gg}$.  It has a corresponding Lie algebra representation,
$d\delta:\tilde{\Gh}\to\operatorname{End}(\tilde{\Gg})$, defined by
$$
\bigl(d\delta(X)\bigr)(\mu)=\left.\frac{d}{dt}\right|_{t=0}
\delta(tX)(\mu).
$$

 From now on, let us fix a dressing orbit ${\mathcal O}$ in $\tilde{\Gg}$.
Let $\mu\in{\mathcal O}$ be a typical element in the orbit (So ${\mathcal O}
=\delta(\tilde{H})(\mu)$.).  We define the stabilizer subgroup by
$$
R_{\mu}=\bigl\{x\in\tilde{H}:\bigl(\delta(x)\bigr)\mu=\mu\bigr\}
\subseteq\tilde{H}.
$$
The corresponding Lie subalgebra is
$$
\Gr_{\mu}=\bigl\{X\in\tilde{\Gh}:\bigl(d\delta(X)\bigr)\mu=0\bigr\}
\subseteq\tilde{\Gh}.
$$
Then the map $\Psi^{\mu}:x\mapsto\delta(x)\mu$ induces a diffeomorphism
$\tilde{H}/{R_{\mu}}\cong{\mathcal O}$.  The range of the differential
map $d\Psi^{\mu}=(d\Psi^{\mu})_e:\tilde{\Gh}\to (T{\mathcal O})_{\mu}$
gives us the tangent space $(T{\mathcal O})_{\mu}$ at $\mu$.  Since
$d\Psi^{\mu}(\tilde{\Gh})=d\delta(\tilde{\Gh})(\mu)$, we have
$\operatorname{Ker}(d\Psi^{\mu})=\Gr_{\mu}$ and $(T{\mathcal O})_{\mu}
={\Gr_{\mu}}^{\perp}$.  And we have the diffeomorphism: $\tilde{\Gh}/
{\Gr_{\mu}}\cong(T{\mathcal O})_{\mu}$.  All this is more or less
the same as the case of the coadjoint orbits (see \cite{CG}).

 Meanwhile, recall that ${\mathcal O}$ is a symplectic leaf, whose
symplectic structure is given by restricting the Poisson bracket on
$\tilde{G}$.  In this way, we obtain a (non-degenerate) skew 2-form
$\omega^{\mu}\in\Lambda^2(T{\mathcal O})_{\mu}$, such that for $X,Y
\in\tilde{\Gh}$, we have:
\begin{align}
&\omega^{\mu}\bigl(d\Psi^{\mu}(X),d\Psi^{\mu}(Y)\bigr)  \notag \\
&=p\cdot(wx'-w'x)+q\cdot(w'y-wy')+\eta_{\lambda}(r)(x\cdot y'-x'\cdot y).
\label{(omega)}
\end{align}
Here $\mu=(p,q,r,s)$ and $X=(x,y,z,w)$, $Y=(x',y',z',w')$.  It is easily
shown to be well-defined (we can show by direct computation).

\begin{lem}
For fixed $\mu$, the above map $\omega^{\mu}\in\Lambda^2(T{\mathcal O})_{\mu}$
is well-defined.
\end{lem}

\begin{proof}
Suppose $X'\in\tilde{\Gh}$ is an arbitrary element such that $d\Psi^{\mu}(X')
=d\Psi^{\mu}(X)$.  Then we can write $X'=X+(a,b,c,d)$, for $(a,b,c,d)\in\Gr_{\mu}$.
By remembering the definition of $\delta$ (from Proposition \ref{dressingGG}),
we can compute the following:
$$
\bigl(d\delta(a,b,c,d)\bigr)(\mu)=\bigl(-dp+\eta_{\lambda}(r)b,
dq-\eta_{\lambda}(r)a,0,p\cdot a-q\cdot b\bigr),
$$
where $\mu=(p,q,r,s)$.  So $(a,b,c,d)\in\Gr_{\mu}$ is characterized by
$$
dp=\eta_{\lambda}(r)b,\quad dq=\eta_{\lambda}(r)a,\quad p\cdot a-q\cdot b=0.
$$
From this and from the definition of $\omega^{\mu}$ given by equation \eqref
{(omega)}, we see easily that for any $Y\in\tilde{\Gh}$, we have
$$
\omega^{\mu}\bigl(d\Psi^{\mu}(X'),d\Psi^{\mu}(Y)\bigr)=\omega^{\mu}
\bigl(d\Psi^{\mu}(X),d\Psi^{\mu}(Y)\bigr).
$$
Similar argument holds for the second entry, and we conclude that $\omega^{\mu}$
is well-defined.
\end{proof}

 We can define $\omega^{\mu}$ for each $\mu\in{\mathcal O}$.  But we can
show that $\omega:\mu\mapsto\omega^{\mu}$ is $\delta(\tilde{H})$-invariant.
To illustrate this more clearly, consider an arbitrary element $\nu$ of
${\mathcal O}$.  Suppose it is written as $\nu=\delta(h)\mu$, for some
$h\in\tilde{H}$.  The following results are true.

\begin{lem}
Let the notation be as above.  Define $\alpha_h:\tilde{H}\to\tilde{H}$
by $\alpha_h(x)=hxh^{-1}$.  Then we have:
$$
\delta(h)\circ\Psi^{\mu}=\Psi^{\nu}\circ\alpha_h.
$$
\end{lem}

\begin{proof}
For any $x\in\tilde{H}$,
\begin{alignat}{2}
\Psi^{\nu}\bigl(\alpha_h(x)\bigr)&=\Psi^{\nu}(hxh^{-1})
=\delta(hxh^{-1})(\nu) &\qquad &({\text {by definition}}) \notag \\
&=\delta(hx)\delta(h^{-1})(\nu) & &(\delta\ {\text {is a left action}})
\notag \\
&=\delta(hx)(\mu) & &(\nu=\delta(h)\mu\ {\text {or}}\ \delta(h^{-1})
\nu=\mu)
\notag \\
&=\delta(h)\delta(x)(\mu)=\delta(h)\Psi^{\mu}(x). & &  \notag
\end{alignat}
\end{proof}

\begin{prop}
Let $\phi_h:=d\bigl(\delta(h)\bigr)_{\mu}:(T{\mathcal O})_{\mu}\to
(T{\mathcal O})_{\nu}$.  We then have, for $X,Y\in\tilde{\Gh}$,
$$
\omega^{\nu}\bigl(\phi_h\circ d\Psi^{\mu}(X),\phi_h\circ d\Psi^{\mu}(Y)
\bigr)=\omega^{\mu}\bigl(d\Psi^{\mu}(X),d\Psi^{\mu}(Y)\bigr).
$$
This illustrates that $\omega:\nu\mapsto\omega^{\nu}$ is
$\delta(\tilde{H})$-invariant.
\end{prop}

\begin{proof}
From the previous lemma, we know that $\delta(h)\circ\Psi^{\mu}=\Psi^{\nu}
\circ\alpha_h$.  By taking differentials, it follows that
$$
\phi_h\circ d\Psi^{\mu}=d\Psi^{\nu}\circ d(\alpha_h).
$$
For convenience, we wrote $d(\alpha_h)=\bigl(d(\alpha_h)\bigr)_e$.
From this, we have
$$
\omega^{\nu}\bigl(\phi_h\circ d\Psi^{\mu}(X),\phi_h\circ d\Psi^{\mu}(Y)
\bigr)=\omega^{\nu}\bigl(d\Psi^{\nu}\circ d(\alpha_h)(X),d\Psi^{\nu}
\circ d(\alpha_h)(Y)\bigr).
$$

To see if it is same as $\omega^{\mu}\bigl(d\Psi^{\mu}(X),d\Psi^{\mu}(Y)
\bigr)$, we will use direct computation.  So let $h=(a,b,c,d)$ and compute
$\alpha_h$, using the multiplication law on $\tilde{H}$ as defined earlier
(in \S1.1):
\begin{align}
&\alpha_h(x,y,z,w)=(a,b,c,d)(x,y,z,w)(a,b,c,d)^{-1} \notag \\
&=(a+e^d x,b+e^{-d}y,c+z+e^{-d}a\cdot y,d+w)(-e^{-d}a,-e^d b,-c+a\cdot b,-d)
\notag \\
&=(a+e^{d}x-e^{w}a,b+e^{-d}y-e^{-w}b,z+a\cdot b+e^{-d}a\cdot y
-e^{d-w}x\cdot b-e^{-w}a\cdot b,w). \notag
\end{align}
It follows that for $X=(x,y,z,w)\in\tilde{\Gh}$, we have:
$$
\bigl(d(\alpha_h)\bigr)(X)=(e^{d}x-wa,e^{-d}y+wb,
z+e^{-d}a\cdot y-e^d x\cdot b+wa\cdot b,w).
$$
Meanwhile, let us write $\mu=(p,q,r,s)$.  Then by definition,
\begin{align}
\nu&=\delta(h)\mu=\bigl(\delta(a,b,c,d)\bigr)(p,q,r,s)  \notag \\
&=\bigl(e^{-d}p+\eta_{\lambda}(r)b,e^{d}q-\eta_{\lambda}(r)a,r,
s+e^{-d}p\cdot a-e^{d}q\cdot b+\eta_{\lambda}(r)a\cdot b\bigr). 
\notag
\end{align}
Therefore, by using the definition of $\omega$ as given in \eqref{(omega)}
and by direct computation using the expressions we obtained above, we see
that for $X=(x,y,z,w)$ and $Y=(x',y',z',w')$, 
\begin{align}
&\omega^{\nu}\bigl(d\Psi^{\nu}\bigl(d(\alpha_h)(X)\bigr),d\Psi^{\nu}
\bigl(d(\alpha_h)(Y)\bigr)\bigr)   \notag \\
&=\bigl(e^{-d}p+\eta_{\lambda}(r)b\bigr)\cdot
\bigl[w(e^{d}x'-w'a)-w'(e^{d}x-wa)\bigr]  \notag \\
&\quad +\bigl(e^{d}q-\eta_{\lambda}(r)a\bigr)\cdot
\bigl[w'(e^{-d}y+wb)-w(e^{-d}y'+w'b)\bigr]  \notag \\
&\quad +\eta_{\lambda}(r)\bigl[(e^{d}x-wa)\cdot(e^{-d}y'+w'b)
-(e^{d}x'-w'a)\cdot(e^{-d}y+wb)\bigr]  \notag \\
&=p\cdot(wx'-w'x)+q\cdot(w'y-wy')+\eta_{\lambda}(r)(x\cdot y'-x'\cdot y)
\notag \\
&=\omega^{\mu}\bigl(d\Psi^{\mu}(X),d\Psi^{\mu}(Y)\bigr).  \notag
\end{align}
Together with the result obtained in the first part, we can conclude our
proof.
\end{proof}

 The invariance of $\omega$ means that it is a $C^{\infty}$ 2-form on
${\mathcal O}$, in its unique symplectic manifold structure inherited
from the Poisson bracket on $\tilde{G}$.

 Our plan (to be carried out in the next section) is to construct a
deformed product on a dense subspace of $C^{\infty}({\mathcal O})$,
by realizing the functions as operators on a Hilbert space.  The
remainder of this section is to make preparations at the level of
orbits.  Let us begin by pointing out that we can identify $\tilde
{H}/{R_{\mu}}$ with the vector space $V=\tilde{\Gh}/{\Gr_{\mu}}$.
What it means is that we are regarding ${\mathcal O}=(T{\mathcal O})_{\mu}$.
See the following remark.

\begin{rem}
In general, we do not know if $R_{\mu}\cong\Gr_{\mu}$ (we may have
$R_{\mu}$ not connected).  But in our case, if we just follow the
definition, it is not difficult to see that they can be identified.
So we will regard $R_{\mu}=\Gr_{\mu}$, for every $\mu$.  Together with
our earlier (spatial) identification $\tilde{H}=\tilde{\Gh}$, we thus
see that $\tilde{H}/{R_{\mu}}=\tilde{\Gh}/{\Gr_{\mu}}$  as vector
spaces.  The advantage of having the description of ${\mathcal O}
(\cong\tilde{H}/{R_{\mu}})$ as a vector space is clear.  We can use
various linear algebraic tools, as well as Fourier transforms.
\end{rem}

 Let us focus our attention to the vector space $V=\tilde{H}/{R_{\mu}}
=\tilde{\Gh}/{\Gr_{\mu}}$.  By the diffeomorphism $\tilde{\Gh}/{\Gr_{\mu}}
\cong(T{\mathcal O})_{\mu}$, it is equipped with the skew, bilinear
form $B^{\mu}$, defined by
$$
B^{\mu}(\dot{X},\dot{Y}):=\omega^{\mu}\bigl(d\Psi^{\mu}(X),d\Psi^{\mu}(Y)
\bigr),
$$
where $X$ and $Y$ are the representatives in $\tilde{\Gh}$ of the classes
$\dot{X},\dot{Y}\in V$.  By the following lemma, we can thus construct
a unique ``self-dual'' measure on $V$.

\begin{lem}\label{selfdualmeasure}
\begin{enumerate}
\item Let $V$ be a (real) vector space, and let $dm$ be a Euclidean
measure on $V$.  Suppose $B(\ ,\ )$ is a nondegenerate bilinear form
on $V$.  Then we may use $B$ to identify $V$ with its dual space $V^*$,
and can define a Fourier transform ${\mathcal F}_B$:
$$
({\mathcal F}_B f)(v)=\int_V f(v')e^{-2\pi iB(v,v')}\,dm(v').
$$
We say that $dm$ is ``self-dual'' if $\|{\mathcal F}_B f\|_2
=\|f\|_2$, where $\|f\|_2^2=\int|f(v)|^2\,dm(v)$.
\item Suppose $\{e_1,e_2,\dots,e_k\}$ is a basis of $V$ and let $dv$
be the normalized Lebesgue measure on $V$.  If we denote by the
same letter $B$ the matrix such that $b_{ij}=B(e_i,e_j)$, then the
self-dual measure on $V$ determined by $B$ is: $dm=|\operatorname
{det}B|^{\frac{1}{2}}\,dv$.
\end{enumerate}
\end{lem}
\begin{proof}
Linear algebra.
\end{proof}

 In our case, on $V$, we have: $dm=|\operatorname{det}(B^{\mu})|^{\frac
{1}{2}}\,d\dot{X}$, where $d\dot{X}$ is the measure inherited from the
Lebesgue measure, $dxdydzdw$, on $\tilde{\Gh}$.  Meanwhile, recall that
$(T{\mathcal O})_{\mu}={{\Gr}_{\mu}}^{\perp}$.  So we may regard
$(T{\mathcal O})_{\mu}$ as the dual vector space of $V=\tilde{\Gh}/
{{\Gr}_{\mu}}$.  Let us give it the dual measure, $d\theta$, of $dm$.
Since $(T{\mathcal O})_{\mu}$ is considered as a subspace ${{\Gr}
_{\mu}}^{\perp}$ of $\tilde{\Gg}$, there already exists a natural
measure on it, denoted by $dl$, inherited from the Plancherel Lebesgue
measure, $dpdqdrds$, on $\tilde{\Gg} (=\tilde{\Gh}^*)$. We will have:
$d\theta=|\operatorname {det} (B^{\mu})|^{-\frac{1}{2}}\,dl$.

 Since $dm$ and $d\theta$ have been chosen by using the (unique)
symplectic structure on ${\mathcal O}$, we know that they will be
the canonical measures on $V(=\tilde{\Gh}/{{\Gr}_{\mu}})$ and
${\mathcal O}\bigl(=(T{\mathcal O})_{\mu}\bigr)$, respectively.
In terms of these canonical measures, there exists the ``symplectic
Fourier transform'', ${\mathcal F}_{\omega}$, from $S_c(V)$ to
$S_c ({\mathcal O})$, as well as its inverse.  That is,
$$
({\mathcal F}_{\omega}f)(l)=\int_{\tilde{\Gh}/{{\Gr}_{\mu}}}
f(\dot{X})\bar{e}\bigl[\langle l,\dot{X}\rangle\bigr]\,dm(\dot{X}),
\qquad f\in S_c(V)
$$
and
$$
({\mathcal F}_{\omega}^{-1}f)(\dot{X})=\int_{(T{\mathcal O})_{\mu}}
f(l)e\bigl[\langle l,\dot{X}\rangle\bigr]\,d\theta(l),\qquad f\in
S_c({\mathcal O}).
$$
The notation $S_c(V)$ means the space of Schwartz functions on $V$
having compact support.  So $S_c(V)\subseteq C_c^{\infty}(V)$.
We can actually work in spaces which are a little larger than
$S_c(V)$ and $S_c({\mathcal O})$, but they are good enough for
our present purposes.

 Let us consider the Hilbert space $L^2(V,dm)$, which contains
$S_c(V)$ as a dense subspace.  By the symplectic Fourier transform,
we see that $L^2(V,dm)\cong L^2({\mathcal O}, d\theta)$.  Due to
the canonical nature of our construction, it is rather easy to see
that the definition of $L^2({\mathcal O},d\theta)$ obtained in this
way does not really depend on the choice of the representative
$\mu\in{\mathcal O}$.  Meanwhile, even though we did not explicitly
mention the case of orbits in ${\mathcal O}(A)$, it is obvious
that everything we have been discussing in this section can be
carried out in exactly the same way (All we need to do is to let
$w$ (or $d$) and $s$ variables to be zero.).

\section{Deformation of the orbits}

 Since a typical dressing orbit ${\mathcal O}$ is a symplectic manifold,
one hopes that there would be a way to define a deformed product on
$C^{\infty}({\mathcal O})$, in the spirit of Weyl quantization and
Moyal products.  Usually, this kind of deformation quantization is
done in terms of ${}^*$-products, involving formal power series
\cite{BFFL}, \cite{Vy}.  Indeed, Arnal and Cortet in \cite{Ar},
\cite{AC} have shown that for nilpotent or some exponential solvable
Lie groups, there exist such quantizations on coadjoint orbits (again
via ${}^*$-products).

 Here, we wish to achieve a similar goal of defining a Moyal-type
deformed product, but without resorting to formal power series.  We
will define our deformed product on a dense subspace $S_c({\mathcal O})$
of $C^{\infty}({\mathcal O})$, by using an operator algebra realization
on a Hilbert space.  The strategy is to relate the deformation with
the irreducible ${}^*$-representation corresponding to the given orbit.
In this article, we plan to discuss only the case of ${\mathcal O}(A)$
in relation with $\operatorname{Irr}(A)$.  The case of ${\mathcal O}
(\tilde{A})$ and $\operatorname{Irr}(\tilde{A})$ will be our future
project (But see Appendix at the end of this section.).

 Let us begin our discussion by turning our attention to the
irreducible ${}^*$-representations.  The next proposition is
very crucial.

\begin{prop}\label{trace}
Let $\pi\in\operatorname{Irr}(A)$ be an irreducible representation
of $A$ acting on ${\mathcal H}_{\pi}$ (See Proposition \ref{repA}
for classification.).  Let ${\mathcal O}\in{\mathcal O}(A)$ be the
corresponding orbit.  Then for $f\in{\mathcal A}$, the operator
$\pi(f)\in{\mathcal B}({\mathcal H}_{\pi})$ turns out to be a
trace-class operator.  Furthermore,
\begin{equation}\label{(trace)}
\operatorname{Tr}\bigl(\pi(f)\bigr)=\int_{\mathcal O}\hat{f}|_{\mathcal O}
(\nu)\,d\theta(\nu),
\end{equation}
where $d\theta$ is the canonical measure on ${\mathcal O}$ defined earlier.
And $\hat{f}$ denotes the partial Fourier transform of $f$ defined by
$$
\hat{f}(p',q',r')=\int f(\tilde{x},\tilde{y},r')\bar{e}(p'\cdot\tilde{x}
+q'\cdot\tilde{y})\,d\tilde{x}d\tilde{y}.
$$
\end{prop}

\begin{rem}
The partial Fourier transform $f\mapsto \hat{f}$ is the usual one (with
respect to the Lebesgue measures), not to be confused with the symplectic
Fourier transform appeared in the previous section.  This can be done
without trouble, since we are (spatially) identifying ${\mathcal O}\bigl
(=(T{\mathcal O})_{\mu}\bigr)$ with the subspace ${\Gr_{\mu}}^{\perp}$
of $\Gg$.
\end{rem}

\begin{proof}
(For $\pi_{p,q}$): This is a trivial case.  Since
$$
\pi_{p,q}(f)=\int f(x,y,0)\bar{e}(p\cdot x+q\cdot y)\,dxdy
=\hat{f}|_{{\mathcal O}_{p,q}}(p,q,0)\in\mathbb{C},
$$
we have: $\operatorname{Tr}\bigl(\pi_{p,q}(f)\bigr)=\hat{f}|_{{\mathcal O}
_{p,q}}(p,q,0)$.

(For $\pi_r$): Note that by using Fourier inversion theorem,
\begin{align}
\bigl(\pi_r(f)\bigr)\xi(u)&=\int f(x,y,r)\bar{e}\bigl[\eta_{\lambda}(r)u
\cdot y\bigr]\xi(u+x)\,dxdy  \notag \\
&=\int\hat{f}|_{{\mathcal O}_r}\bigl(\tilde{p},\eta_{\lambda}(r)u,r\bigr)
\xi(u+\tilde{x})e(\tilde{p}\cdot\tilde{x})\,d\tilde{p}d\tilde{x}
\notag \\
&=\int K(u,\tilde{x})\xi(\tilde{x})\,d\tilde{x},  \notag
\end{align}
where $K(u,\tilde{x})=\int\hat{f}|_{{\mathcal O}_r}\bigl(\tilde{p},
\eta_{\lambda}(r)u,r\bigr)e\bigl[\tilde{p}\cdot(\tilde{x}-u)
\bigr]\,d\tilde{p}$.  That is, $\pi_r(f)$ is an integral operator
whose kernel is given by $K$, which is clearly an $L^2$-function
since $\hat{f}$ is a Schwartz function.  This means that $\pi_r(f)$
is a trace-class operator.  Moreover,
\begin{align}
\operatorname{Tr}\bigl(\pi_r(f)\bigr)&=\int K(u,u)\,du=\int
\hat{f}|_{{\mathcal O}_r}\bigl(\tilde{p},\eta_{\lambda}(r)u,
r\bigr)\,d\tilde{p}du   \notag \\
&=\int\hat{f}|_{{\mathcal O}_r}(\tilde{p},u,r)\frac{1}{\bigl|\eta_
{\lambda}(r)\bigr|^n}\,d\tilde{p}du.  \notag
\end{align}
But the measure $\bigl|\eta_{\lambda}(r)\bigr|^{-n}\,d\tilde{p}du$
is none other than the canonical measure $d\theta(\tilde{p},u,r)$ on
${\mathcal O}_r$.  To see this more clearly, recall the definition
of the bilinear map $B^{\mu}$ on $\Gh/{{\Gr}_{\mu}}\cong{\mathcal O}$.
In our case, we may choose $\mu=(0,0,r)$ and $\Gh/{{\Gr}_{\mu}}=\Gh/\Gz$.
We then have:
$$
B^{\mu}(\dot{X},\dot{Y})=\eta_{\lambda}(r)(x\cdot y'-x'\cdot y),
$$
where $\dot{X}$ and $\dot{Y}$ are the classes in $\Gh/\Gz$ represented
by $X=(x,y,0)$ and $Y=(x',y',0)$.  By simple calculation, we have:
$\operatorname{det}(B^{\mu})=\bigl(\eta_{\lambda}(r)\bigr)^{2n}$.  From
this, it follows from the discussion following Lemma \ref{selfdualmeasure}
that the canonical measure on ${\mathcal O}_r$ is:
$$
d\theta(p',q',r)=\bigl|\operatorname{det}(B^{\mu})\bigr|^{-\frac{1}{2}}
\,dp'dq'=\bigl|\eta_{\lambda}(r)\bigr|^{-n}\,dp'dq'.
$$
This verifies the trace formula: equation \eqref{(trace)}.
\end{proof}

\begin{cor}
Let $\pi\in\operatorname{Irr}(A)$.  Since each $\pi(f)$, $f\in
{\mathcal A}$, is a trace-class operator, it is also Hilbert--Schmidt.
In our case, the Hilbert--Schmidt norms are given by
\begin{align}
&\bigl\|\pi_{p,q}(f)\bigr\|_{\operatorname{HS}}^2=\bigl|\hat{f}
(p,q,0)\bigr|^2,  \notag \\
&\bigl\|\pi_r(f)\bigr\|_{\operatorname{HS}}^2=\int\overline{\hat{f}
(p',q',r)}f(p',q',r)\bigl|\eta_{\lambda}(r)\bigr|^{-n}\,dp'dq'
=\int\bigl|\hat{f}(p',q',r)\bigr|^2\,d\theta.  \notag
\end{align}
Indeed, we actually have, for $f,g\in{\mathcal A}$:
$$
\operatorname{Tr}\bigl(\pi(g)^*\pi(f)\bigr)=\operatorname{Tr}\bigl
(\pi(g^*\times f)\bigr)=\bigl\langle\hat{f}|_{{\mathcal O}_r},
\hat{g}|_{{\mathcal O}_r}\bigr\rangle_{\mathcal O},
$$
where $\times$ is the multiplication on ${\mathcal A}$ and $g\mapsto
g^*$ is the involution on ${\mathcal A}$, while $\langle\ ,\ \rangle_
{\mathcal O}$ denotes the inner product on $L^2({\mathcal O},d\theta)$.
\end{cor}

\begin{proof}
We just need to remember the definitions of the multiplication and
involution on ${\mathcal A}$ (for instance, see Propositions 2.8 and 2.9
of \cite{BJKp2}), and use the trace formula obtained in the previous
proposition.  The result follows from straightforward computation.
\end{proof}

 Since the elements $\pi(f)$, $f\in{\mathcal A}$, form a dense subspace
in $HS({\mathcal H}_{\pi})$, the above Corollary implies that we have a
Hilbert space isomorphism between $HS({\mathcal H}_{\pi})$ and $L^2
({\mathcal O},d\theta)$.  To be a little more precise, let us consider
the map from $HS({\mathcal H}_{\pi})$ to $L^2({\mathcal O},d\theta)$
by naturally extending the map $\pi(f)\mapsto\hat{f}|_{\mathcal O}$.
By the result we just obtained, the map is an isometry, preserving
the inner product.  It is clearly onto.  Let us from now on consider
its inverse map and denote it by $S_{\pi}$.  In this way, we have
the spatial isomorphism, $S_{\pi}:L^2({\mathcal O},d\theta)\cong HS
({\mathcal H}_{\pi})$.

 We are now ready to discuss the deformation of the orbits ${\mathcal O}$.
The point is that through the map $S_{\pi}$, an arbitrary element $\phi\in
S_c({\mathcal O})$ can be considered as a (Hilbert--Schmidt) operator on
${\mathcal H}_{\pi}$.  Let us denote this correspondence by $Q_{\pi}$.  It
is the same as the map $S_{\pi}$ above, but now we work with the operator
norm on ${\mathcal B}({\mathcal H}_{\pi})$ instead of the Hilbert--Schmidt
norm on $HS({\mathcal H}_{\pi})$.  The result is summarized in the below.

\begin{prop}\label{moyalproductA}
Let $Q_{\pi}:S_c({\mathcal O})\to{\mathcal B}({\mathcal H}_{\pi})$ be
defined as above.  Then $Q_{\pi}$ determines a deformed product on $S_c
({\mathcal O})$ by $Q_{\pi}(\phi\times_Q\psi)=Q_{\pi}(\phi)Q_{\pi}(\psi)$.
The involution on $S_c({\mathcal O})$ will be given by $Q_{\pi}(\psi^*)
={Q_{\pi}(\psi)}^*$.  They are described in the below:
\begin{enumerate}
\item (The case of the 1-point orbit ${\mathcal O}_{p,q}$): For $\phi,\psi
\in S_c({\mathcal O}_{p,q})$,
\begin{align}
\bigl(\phi\times_Q\psi\bigr)(p,q,0)&=\phi(p,q,0)\psi(p,q,0),  \notag\\
\psi^*(p,q,0)&=\overline{\psi(p,q,0)}.  \notag
\end{align}
\item (The case of the $2n$-dimensional orbit ${\mathcal O}_r$): For
$\phi,\psi\in S_c({\mathcal O}_{r})$, 
\begin{align}
&\bigl(\phi\times_Q\psi\bigr)(\alpha,\beta,r)
=\int\phi(\tilde{p},\beta,r)\psi\bigl(\alpha,\beta+\eta_{\lambda}(r)
\tilde{x},r\bigr)e\bigl[(\tilde{p}-\alpha)\cdot\tilde{x}\bigr]\,
d\tilde{p}d\tilde{x},  \notag \\
&\psi^*(\alpha,\beta,r)=\int\overline{\psi(\tilde{p},\tilde{q},r)}
\bar{e}\bigl[(\alpha-\tilde{p})\cdot x+(\beta-\tilde{q})\cdot y\bigr]
\bar{e}\bigl[\eta_{\lambda}(r)x\cdot y\bigr]\,d\tilde{p}d\tilde{q}dxdy.
\notag
\end{align}
\end{enumerate}
\end{prop}

\begin{proof}
The computations are straightforward from the definitions.  So we will
just verify the multiplication formula $\phi\times_Q\psi$, for the
case of the $2n$-dimensional orbit.  Indeed, for $\xi\in{\mathcal B}
\bigl(L^2(\mathbb{R}^n)\bigr)$, we have:
\begin{align}
&Q_{{\pi}_r}(\phi)Q_{{\pi}_r}(\psi)\xi(u)  \notag \\
&=\int\phi(\tilde{p},\tilde{q},r)e(\tilde{p}\cdot x+\tilde{q}\cdot y)\bar{e}
\bigl[\eta_{\lambda}(r)u\cdot y\bigr]Q_{{\pi}_r}(\psi)\xi(u+x)\,d\tilde{p}
d\tilde{q}dxdy  \notag \\
&=\int\phi\bigl(\tilde{p},\eta_{\lambda}(r)u,r\bigr)e(\tilde{p}\cdot x)
Q_{{\pi}_r}(\psi)\xi(u+x)\,d\tilde{p}dx  \notag \\
&=\int\phi\bigl(\tilde{p},\eta_{\lambda}(r)u,r\bigr)e(\tilde{p}\cdot x)
\psi\bigl(\tilde{\tilde{p}},\eta_{\lambda}(r)(u+x),r\bigr)e(\tilde{\tilde{p}}
\cdot x')\xi(u+x+x')\,d\tilde{p}dxd\tilde{\tilde{p}}dx'  \notag \\
&=\int(\phi\times_Q\psi)\bigl(\tilde{\tilde{p}},\eta_{\lambda}(r)u,r\bigr)
e(\tilde{\tilde{p}}\cdot x')\xi(u+x')\,d\tilde{\tilde{p}}dx'  \notag \\
&=Q_{{\pi}_r}(\phi\times_Q\psi)\xi(u),  \notag
\end{align}
where
$$
(\phi\times_Q\psi)\bigl(\tilde{\tilde{p}},\eta_{\lambda}(r)u,r\bigr)
=\int\phi\bigl(\tilde{p},\eta_{\lambda}(r)u,r\bigr)\psi\bigl(\tilde
{\tilde{p}},\eta_{\lambda}(r)(u+x),r\bigr)e\bigl[(\tilde{\tilde{p}}
-\tilde{p})\cdot x\bigr]\,d\tilde{p}dx.
$$
\end{proof}

\begin{rem}
Since the ${}^*$-algebra structure on $S_c({\mathcal O})$ has been defined
via a ${}^*$-representation, the properties like associativity of the
multiplication are immediate.  Although we do not plan to point out the
actual deformation process, the product on $S_c({\mathcal O})$ as obtained
above can be shown to be a deformation quantization of the pointwise
product on $S_c({\mathcal O})$, in the direction of the symplectic
structure (given by $\omega^{\mu}$ or $B^{\mu}$ in our case) on the orbit
${\mathcal O}$: For instance, just as in \cite{Rf3} or in \cite{BJKp1},
we can replace $B^{\mu}(\dot{X},\dot{Y})$ by $\frac{1}{\hbar}B^{\mu}(\hbar
\dot{X},\hbar\dot{Y})$ and proceed, with $\hbar$ being the deformation
parameter.
\end{rem}

 In this sense, we call $\phi\times_Q\psi$ a Moyal-type product, because
it resembles the process of Weyl quantization of $C^{\infty}(\mathbb
{R}^{2k})$ and Moyal product.  Recall that in Weyl quantization
(e.\,g. see \cite{Fo}), functions $\phi,\psi\in C_c^{\infty}(\mathbb
{R}^{2k})$ are associated to certain operators $F_{\phi}$ and $F_{\psi}$
involving Schr\"odinger's $P$, $Q$ operators, and the operator
multiplication $F_{\phi}F_{\psi}$ is defined to be $F_{\phi\times\psi}$,
giving us the Moyal product.  Several authors since have modified this
process to obtain various versions of Moyal-type products (mostly in
terms of ${}^*$-products and formal power series).  Above formulation
is just one such.  On the other hand, note that in our case, we do not
have to resort to the formal power series.  Our quantization is carried
out using the $C^*$-algebra framework.

\begin{defn}
We will write $A_{\pi}:=\overline{Q_{\pi}\bigl(S_c({\mathcal O}_{\pi})
\bigr)}^{\|\ \|_{\operatorname{op}}}$, the norm closure in ${\mathcal B}
({\mathcal H}_{\pi})$ of the ${}^*$-algebra $S_c({\mathcal O}_{\pi})$
considered above.  The $C^*$-algebras $A_{\pi}$ will be considered as
the quantizations of the orbits ${\mathcal O}_{\pi}$.
\end{defn}

 In view of the proposition \ref{moyalproductA}, we may as well say that
each irreducible representation of $(A,\Delta)$ ``models'' the deformed
multiplication on each $S_c({\mathcal O})$.  This is not necessarily a
very rigorous statement, but it does give us a helpful insight: Note
that in the geometric quantization program, especially in the program
introduced in \cite{BFFL}, one studies the representation theory of Lie
groups (or more general objects) via deformed products of functions.

 Meanwhile, from the Corollary to Proposition \ref{trace}, the following
Plancherel-type result is immediate.  Here, ${\mathcal A}$ is viewed as
a dense subspace of ${\mathcal H}$, which is the Hilbert space consisting
of the $L^2$-functions in the $(x,y,r)$ variables.  It is the Hilbert space
on which our Hopf $C^*$-algebra $(A,\Delta)$ acts (See \S1.2, as well as
our previous paper \cite{BJKp2}.).

\begin{prop}
For $f\in{\mathcal A}$, viewed as a Schwartz function contained in
the Hilbert space ${\mathcal H}$, we have:
\begin{equation}\label{(plancherel)}
\|f\|_2^2=\int\bigl\|\pi_r(f)\bigr\|_{\operatorname{HS}}^2
\bigl|\eta_{\lambda}(r)\bigr|^n\,dr.
\end{equation}
\end{prop}

\begin{proof}
For $f\in{\mathcal A}$,
\begin{align}
\|f\|_2^2=\|\hat{f}\|_2^2&=\int\bigl|\hat{f}(p,q,r)\bigr|^2\,dpdqdr
\notag \\
&=\int\left(\int\bigl|\hat{f}(p,q,r)\bigr|^2\bigl|\eta_{\lambda}
(r)\bigr|^{-n}\,dpdq\right)\bigl|\eta_{\lambda}(r)\bigr|^n\,dr
\notag \\
&=\int\bigl\|\pi_r(f)\bigr\|_{\operatorname{HS}}^2\bigl|
\eta_{\lambda}(r)\bigr|^n\,dr.   \notag
\end{align}
\end{proof}

 In $\operatorname{Irr}(A)\cong{\mathcal O}(A)$, the topology on
the subset $\{\pi_r\}$ (consisting of generic representations)
is the standard one, which is just the vector space topology on
a line (the $r$-axis) with one point removed (at $r=0$).  And,
the $\{\pi_r\}$ form a dense subset in $\operatorname{Irr}(A)$.
So $d\mu=\bigl|\eta_{\lambda}(r)\bigr|^n\,dr$ in \eqref{(plancherel)}
may be regarded as a measure on $\operatorname {Irr}(A)$, which is
(densely) supported on $\{\pi_r\}\subseteq\operatorname{Irr}(A)$.  
Keeping the analogy with the representation theory of nilpotent Lie
groups \cite{CG}, we will call $d\mu$ the {\em Plancherel measure\/}
on $\operatorname{Irr}(A)$.

 Relative to the Plancherel measure $d\mu$, we can construct the
following direct integral of Hilbert spaces:
$$
\int_{\operatorname{Irr}(A)}^{\oplus}\operatorname{HS}
({\mathcal H}_{\pi})\,d\mu(\pi),
$$
where $\operatorname{HS}({\mathcal H}_{\pi})$ denotes the space of
Hilbert--Schmidt operators on ${\mathcal H}_{\pi}$.  The Plancherel
formula \eqref{(plancherel)} implies that the map
$$
T:{\mathcal A}
\to\int_{\operatorname{Irr}(A)}^{\oplus}\operatorname{HS}
({\mathcal H}_{\pi})\,d\mu(\pi)
$$
defined by $T(f)(\pi):=\pi(f)$, $\pi\in\operatorname{Irr}(A)$,
is an isometry.  Clearly, $T$ will extend to an isometry on
${\mathcal H}$.  Actually, $T$ is an onto isometry.

\begin{theorem}\label{Plancherel}(The Plancherel Theorem)

\noindent The map $T:{\mathcal H}\to\int_{\operatorname{Irr}(A)}
^{\oplus}\operatorname{HS}({\mathcal H}_{\pi})\,d\mu(\pi)$ as
defined above is an onto isometry.  In other words, the spatial
isomorphism $T$ gives a direct integral decomposition of
${\mathcal H}$, via the Plancherel formula \eqref{(plancherel)}.
\end{theorem}

\begin{proof}
This is a version of the Plancherel Theorem.  The proof can be given
following the direct integral analysis of Dixmier \cite[\S18]{D1}.
We can further say that the Plancherel measure is actually unique.
This result illustrates the fact that $A$ is a type I $C^*$-algebra.
\end{proof}

 Recall that on ${\mathcal H}$, the algebra ${\mathcal A}$ (or $A$)
acts by regular representation $L$.  The precise definition can
be found in our previous papers, but $L(f)$, $f\in{\mathcal A}$, is
none other than the multiplication operator defined by $L(f)\xi=
f\times\xi$, where $\xi\in{\mathcal A}\subseteq{\mathcal H}$ and
$\times$ is the multiplication on ${\mathcal A}$.  Meanwhile, on
each fiber $\operatorname{HS}({\mathcal H}_{\pi})$ of the decomposition,
$f\in{\mathcal A}$ acts by $F\mapsto\pi(f)F$, where the right hand side
means the operator multiplication between the two Hilbert--Schmidt
operators on ${\mathcal H}_{\pi}$.

 It is easy to see that $T$ intertwines these actions.  That is,
if $f,\xi\in{\mathcal A}(\subseteq{\mathcal H})$, then
$$
T\bigl(L(f)\xi\bigr)(\pi)=T(f\times\xi)(\pi)=\pi(f\times\xi)
=\pi(f)\pi(\xi)=\pi(f)\bigl(T(\xi)(\pi)\bigr),
$$
which holds for all $\pi$.  In other words, we have the equivalence
of representations:
$$
L\cong\int_{\operatorname{Irr}(A)}^{\oplus}\pi\,d\mu(\pi).
$$
This means that the regular representation has a direct integral
decomposition into irreducible representations.

\begin{rem}
Even though the Plancherel measure is supported only on the set of
generic irreducible representations $\{\pi_r\}$, note that since $\mu
\bigl(\operatorname{Irr}(A)\setminus\{\pi_r\}\bigr)=0$ and closure
of the $\{\pi_r\}$ is all of $\operatorname{Irr}(A)$, the above result
is consistent with the observation made in our earlier paper (in \cite
{BJKp2}) on the amenability of $A$.  Meanwhile, we point out here that
very recently, Desmedt \cite{Dem} has studied the amenability problems
and a generalized version of Plancherel theorem, in the setting of locally
compact quantum groups (The author thanks the referee for this information.).
\end{rem}

 Before we wrap up, let us point out the following interesting
observation: By our quantization map $Q_{\pi}$, we saw that
the deformed product on $S_c({\mathcal O})$ is actually an operator
multiplication, and we summarized this situation by saying that
``each irreducible representation models the deformed multiplication
on each $S_c({\mathcal O})$''.  Now by the Plancherel theorem, the
regular representation $L$ of $A$ has a direct integral decomposition
into the irreducible representations.  Since the regular representation
is given by the left multiplication on ${\mathcal A}$ (and $A$), and
since each irreducible representation models the Moyal-type deformed
multiplication on each $S_c({\mathcal O})$, the Plancherel theorem can
be loosely stated as follows: ``the twisted product on $A$ is patched-up
from the deformed multiplications on the $S_c({\mathcal O})$''.

 In our case, the twisted product on $A$ was obtained directly as
a deformation quantization of the Poisson bracket, via a certain
(continuous) cocycle \cite{BJKp1}, \cite{BJKp2}.  The above paragraph
suggests a more geometric approach such that one may try to construct
the twisted product by first studying the individual dressing orbits
(symplectic leaves), find a Moyal-type products on them, and then
``patch-up'' these deformed products.  A similar idea is being used
in \cite{NatNP}, although the settings are different from ours.

 In general, we do not expect it to work fully, due to various
obstructions caused by the complexities of the symplectic leaves
themselves and of the way the leaves lie inside the Poisson manifold.
Nevertheless, this observation underlines the point that the dressing
orbits and the representation theoretical analysis play a very useful
role in the development of quantization methods.

\section{Appendix: Deformation of the orbits in ${\mathcal O}(\tilde{A})$}

 Finally, a short remark is in order for the case of the orbits
in ${\mathcal O}(\tilde{A})$, in relation with the representations
in $\operatorname{Irr}(\tilde{A})$.  Generally speaking, the ideas
we followed in section 4 do go through, in the sense that we can
define a deformed multiplication on an orbit ${\mathcal O}\in
{\mathcal O}(\tilde{A})$, which is modeled by the irreducible
representation corresponding to the orbit.  In addition, the
Plancherel type result exists, giving us an interpretation as
above that the regular representation is ``patched-up'' of the
deformed multiplications on the orbits.

 On the other hand, not all the steps go through and some modifications
should be made.  We do not plan to give any detailed discussion here
(which would be rather lengthy and since it is still in the works),
but we will briefly indicate in the below where the modifications should
occur.

 For the cases of the representations $\tilde{\pi}_s$ and $\tilde
{\pi}_{p,q}$, essentially the same results hold as in Propositions
\ref{trace} and \ref{moyalproductA}.  That is, we do have the spatial
isomorphism $L^2({\mathcal O}_{\pi},d\theta)\cong\operatorname{HS}
({\mathcal H}_{\pi})$, and from this the Moyal-type deformed multiplication
on the orbits can be obtained.  The case of the representations $\tilde
{\pi}_{r,s}$ is when we need some care: The operators $\tilde{\pi}_{r,s}
(f)$, $f\in\tilde{\mathcal A}$, are no longer trace-class, and it seems
we need to incorporate a kind of a ``formal degree'' operator (as in
\cite{Dur}) for $\tilde{\pi}_{r,s}$, to define the Hilbert-Schmidt
operators.  Even with this adjustment, the Hilbert-Schmidt operator
space is not isomorphic to $L^2({\mathcal O}_{r,s},d\theta)$.

 Some of these are serious obstacles, but it turns out that there is
still a way to define a deformed multiplication on each orbit, again
modeled by the irreducible representations.  We do not seem to have 
spatial isomorphisms between the $L^2({\mathcal O}_{r,s}, d\theta)$
and the $\operatorname{HS}({\mathcal H}_{\pi_{r,s}})$, but we can still
show that $\int^{\oplus}L^2({\mathcal O}_{r,s},d\theta)\,ds\cong
\int^{\oplus}\operatorname{HS}({\mathcal H}_{\pi_{r,s}})\,ds$. From
this, it follows that: $\int^{\oplus}L^2({\mathcal O}_{r,s},d\theta)
\,\bigl|\eta_{\lambda}(r)\bigr|^{-1}drds\cong\tilde{\mathcal H}$, where
$\tilde{\mathcal H}$ is the Hilbert space on which our $C^*$-algebra
$\tilde{A}$ acts by regular representation (as in Example 3.6 of
\cite{BJKp2}).  This would be the result taking the place of Theorem
\ref{Plancherel} above.


\bibliography{refpp}

\providecommand{\bysame}{\leavevmode\hbox to3em{\hrulefill}\thinspace}
\providecommand{\MR}{\relax\ifhmode\unskip\space\fi MR }
\providecommand{\MRhref}[2]{%
  \href{http://www.ams.org/mathscinet-getitem?mr=#1}{#2}
}
\providecommand{\href}[2]{#2}
\begin{thebibliography}{10}

\bibitem{Ar}
D.~Arnal, \emph{${}^*$ products and representations of nilpotent groups},
  Pacific J. Math. \textbf{114} (1984), no.~2, 285--308.

\bibitem{AC}
D.~Arnal and J.~C. Cortet, \emph{Repr\'esentations {${}^*$} des groupes
  exponentiels}, J. Funct. Anal. \textbf{92} (1990), 103--135 (French).

\bibitem{BFFL}
F.~Bayen, M.~Flato, C.~Fronsdal, A.~Lichnerowicz, and D.~Sternheimer,
  \emph{Deformation theory and quantization {I}, {II}}, Ann. Phys. \textbf{110}
  (1978), 61--110, 111--151.

\bibitem{CG}
L.~Corwin and F.~P. Greenleaf, \emph{Representations of {N}ilpotent {L}ie
  {G}roups and {T}heir {A}pplications. {P}art 1}, Cambridge studies in advanced
  mathematics, no.~18, Cambridge Univ. Press, 1990.

\bibitem{Dem}
P.~Desmedt, 2003, Ph.D. thesis (KU Leuven, Belgium).

\bibitem{D1}
J.~Dixmier, \emph{{$C^*$}-algebras}, North-Holland, 1977.

\bibitem{Dur}
M.~Duflo and M.~Ra{\"i}s, \emph{Sur l'analyse harmonique sur les groupes de
  {L}ie r{\'e}solubles}, Ann. Scient. \'Ec. Norm. Sup., $4^e$ s\'erie
  \textbf{t. 9} (1976), 107--144 (French).

\bibitem{Fo}
G.~B. Folland, \emph{Harmonic {A}nalysis in {P}hase {S}pace}, Annales of
  Mathematical Studies, no. 122, Princeton University Press, 1989.

\bibitem{BJKp1}
B.~J. Kahng, \emph{Deformation quantization of certain non-linear {P}oisson
  structures}, Int. J. Math. \textbf{9} (1998), 599--621.

\bibitem{BJKp2}
\bysame, \emph{Non-compact quantum groups arising from {H}eisenberg type {L}ie
  bialgebras}, J. Operator Theory \textbf{44} (2000), 303--334.

\bibitem{BJKhj}
\bysame, \emph{${}^*$-representations of a quantum {H}eisenberg group algebra},
  Houston J. Math. \textbf{28} (2002), 529--552.

\bibitem{BJKppha}
\bysame, \emph{Haar measure on a locally compact quantum group}, J. Ramanujan
  Math. Soc. \textbf{18} (2003), 385--414.

\bibitem{Ki2}
A.~A. Kirillov, \emph{Unitary representations of nilpotent {L}ie groups},
  Russian Math. Surveys \textbf{17} (1962), no.~4, 53--104, Translated from
  Usp. Mat. Nauk. {\bf 17} (1962), 57--110.

\bibitem{Ki}
\bysame, \emph{Elements of the {T}heory of {R}epresentations}, Springer-Verlag,
  Berlin, 1976.

\bibitem{Ki4}
\bysame, \emph{The orbit method, {I}: {G}eometric quantization}, Representation
  {T}heory of {G}roups and {A}lgebras, Contemp. Math., no. 145, American
  Mathematical Society, 1993, pp.~1--32.

\bibitem{Ki5}
\bysame, \emph{Merits and demerits of the orbit method}, Bull. AMS \textbf{36}
  (1999), no.~4, 433--488.

\bibitem{Ko}
B.~Kostant, \emph{Quantization and unitary representations}, Lecture Notes in
  Modern Analysis and Applications III, Lecture Notes in Math., no. 170,
  Springer-Verlag, 1970, pp.~87--208.

\bibitem{KuVa}
J.~Kustermans and S.~Vaes, \emph{Locally compact quantum groups}, Ann. Scient.
  \'Ec. Norm. Sup., $4^e$ s\'erie \textbf{t. 33} (2000), 837--934.

\bibitem{LL}
H.~Leptin and J.~Ludwig, \emph{Unitary {R}epresentation {T}heory of
  {E}xponential {L}ie {G}roups}, Walter de Gruyter \& Co., 1994.

\bibitem{LS}
S.~L. Levendorskii and Y.~S. Soibelman, \emph{Algebras of functions on compact
  quantum groups, {S}chubert cells, and quantum tori}, Comm. Math. Phys.
  \textbf{139} (1991), 141--170.

\bibitem{LW}
J.~H. Lu and A.~Weinstein, \emph{Poisson {L}ie groups, dressing transformations
  and {B}ruhat decompositions}, J. Diff. Geom. \textbf{31} (1990), 501--526.

\bibitem{NatNP}
T.~Natsume, R.~Nest, and I.~Peter, \emph{{$C^*$}-algebraic deformation
  quantization of symplectic manifolds}, 1999, preprint.

\bibitem{Rf3}
M.~A. Rieffel, \emph{Lie group convolution algebras as deformation
  quantizations of linear {P}oisson structures}, Amer. J. Math. \textbf{112}
  (1990), 657--685.

\bibitem{Rf4}
\bysame, \emph{Deformation quantization for actions of {$R^d$}}, Memoirs of the
  AMS, no. 506, American Mathematical Society, Providence, RI, 1993.

\bibitem{Se}
M.~A. Semenov-Tian-Shansky, \emph{Dressing transformations and {P}oisson group
  actions}, Publ. RIMS, Kyoto Univ. \textbf{21} (1985), 1237--1260.

\bibitem{SV}
Y.~S. Soibelman and L.~L. Vaksman, \emph{Algebra of functions on the quantum
  group {$SU(2)$}}, Functional Anal. Appl. \textbf{22} (1988), 170--181.

\bibitem{Vy}
J.~Vey, \emph{D\'eformation du crochet de {P}oisson sur une vari\'et\'e
  symplectique}, Comm. Math. Helv. \textbf{50} (1975), 421--454 (French).

\end{thebibliography}

\bibliographystyle{amsplain}

\end{document}